\input amstex
\documentstyle{amsppt}
\magnification=1200
\overfullrule=0pt
\NoRunningHeads
\nologo

\define\mbn{\overline{M}_{0n}}
\define\mbs{\overline{M}_{0S}}

\define\s{\sigma}
\define\ol{\overline}

\define\on {\overline{n}}
\define\onn {\overline {n-1}}
\define\n {\{1,\ldots ,n \}}
\define\vp {\varphi}
\define\med{\medskip}
\define\sm{\smallskip}

\topmatter
\title
The intersection form in $H^*(\mbn)$ and the explicit 
K\"unneth formula in quantum cohomology
\endtitle


\author 
Ralph Kaufmann 
\endauthor

\affil 
Max--Planck--Institut f\"ur Mathematik, Bonn, Germany 
\endaffil
\abstract
We prove a general formula for the intersection form of 
two arbitrary monomials in boundary
divisors. Furthermore we present a tree basis of the cohomology of $\mbn$. With
the help of the intersection form we determine the Gram matrix for this basis
and give a formula for its inverse. This enables us to calculate the tensor product of the higher order
multiplications arising in quantum cohomology and formal Frobenius manifolds.
In the context of quantum cohomology this establishes the explicit 
K\"unneth formula.
\endabstract
\endtopmatter


\document
\TagsOnRight

{\bf 0. Introduction}

\smallskip

Let $\mbn$ be the moduli space of genus $0$ curves with $n$ marked points.
Its cohomology ring was determined by Keel [Ke], who gave a presentation 
in terms of boundary divisors, their intersections and relations. A boundary
divisor is specified by a 2-partition $S_1 \amalg S_2$ of 
$\on :=\{1,\dots ,n\}$.
The additive structure of this ring was studied and presented in [KM] and [KMK].
Although much about the structure of this ring is known there are still 
several open 
questions. The complete study of the intersection theory of this space 
however is
of importance for the theory of quantum cohomology. In particular it
is necessary in order to understand the K\"unneth formula for 
quantum cohomology,
which is given by the tensor product of Cohomological Field Theories, 
cf\. [KM] and [KMK].

In \S 2 of this paper we prove a formula for the intersection form
for any two polynomials in the boundary divisors of complementary degree.
More precisely, after the introduction of the notion of trees with 
multiplicities and
good multiplicity orientations we can formulate the following 
\proclaim {Theorem}
Let $mon(\s_1,m_1)$ and $mon(\s_2,m_2)$ be two monomials of complementary
degree in $H^*(\mbn)$. If there is no good multiplicity orientation of
$(\tau,m) := \tau(\s_1 \cup \s_2,m_1 +m_2)$ then
$\langle mon(\s_1,m_1) mon(\s_1,m_1) \rangle =0$. If there does exist one then:
$$
\langle mon(\s_1,m_1) mon(\s_1,m_1) \rangle 
= \prod_{v \in V_{\tau}} (-1)^{|v|-3}
\frac {(|v|-3)!} {\prod_{f \in F(v)} (mult (f))!^2}
\prod_{e\in E_{\tau}}(m(e)-1)!,
$$
where $mult$ is the unique multiplicity orientation for $(\tau,m)$
provided by the Lemma 2.3 whose value is given in the formula (2.3).
\endproclaim
 The notation 
$mon ((\tau,m)):= \prod_{e \in E_{\tau}} D_{\s(e)}^{m(e)}$ 
used in this theorem 
along with  an exposition of the different
combinatorics of trees involved in the intersection theory of $\mbn$ is explained in 
the introductory \S 1, where also Keel's presentation is briefly reviewed.

In \S 3 we will give a monomial basis ${\Cal B}_n$ of $H^*(\mbn)$ together with
a tree representation for it. The basis which is presented here and the proof
of linear independence is
inspired by the work of Yuzvinsky [Yu], who worked out a basis in another 
presentation of the cohomology ring developed by DeConcini and Procesi [CP] via 
hyperplane arrangements. Using the results of \S 2 we can write down the Gram 
matrix
for this basis and give a formula for its inverse.

The realization of this basis
in terms of boundary divisors is necessary for applications to quantum 
cohomology and operads [KM,G], since these structures make explicit use of the 
presentation
of $H^*(\mbn)$ in terms of tree strata.

As an application to this field we use the results of \S 3 to calculate the 
tensor product of the higher order products
and correlation functions stemming from a tree level cohomological field theory
which appear in the tensor product of formal Frobenius
manifolds and yield the explicit K\"unneth formula for
quantum cohomology.

\proclaim {Corollary}
For two projective algebraic manifolds $V$ and $W$
the potential $\Phi^{V\times W}$ yielding the quantum cohomology of
$V \times W$ in terms of $\Phi^V$ and $\Phi^W$ is given by the formula:
$$
\Phi^{V\times W}(\gamma' \otimes \gamma'') = \sum_{n\geq 3} \frac {1}{n!} 
\sum_{\mu, \nu \in {\Cal B}_n} 
Y'(\check {\mu}) (\gamma^{\prime \otimes n}) m_{\mu \nu}
Y''(\check{\nu}) (\gamma^{\prime \prime \otimes n})
$$
where ${\Cal B}_n$ is the basis of \S 3, 
$(m_{\mu \nu})_{\mu \nu \in{\Cal B}_n} $ is its Gram matrix (3.14),
$\{\check{\mu}|\mu \in {\Cal B}_n\}$ is the dual basis obtained via the inverse 
Gram matrix (3.15) and  $\{Y'(\tau)\}$ resp.\ $\{Y''(\tau)\}$ are the operadic 
ACF's obtained from $\Phi^V$ resp.\ $\Phi^W$ via (4.3) and (4.5).
\endproclaim

As an example the first higher order products are written out explicitly.

I would like to thank S\. Yuzvinsky for sending his manuscript before 
publication and the Max-Planck-Institut f\"ur Mathematik for its financial 
support and stimulating atmosphere.
Most of all I want to express my gratitude to Yu\. I\. Manin for his continuous 
support and encouragement. 

\medskip

{\bf 0.1 Notation.}

\smallskip

Throughout this paper we will denote by $\subset$ the strict 
inclusion and use $\subseteq$ for the not necessarily strict relation. 
Furthermore we denote
by $\Bbb N$ the positive integers, and we will use the notation $\on$ to denote 
the set $\n$. 
\bigskip

\newpage

{\bf \S 1 Partitions and trees}

\medskip

{\bf 1.1 Notation.}

\smallskip

We will consider a tree $\tau$ to be a collection of sets of
vertices, edges and tails  
$(V_{\tau},E_{\tau}, T_{\tau})$ with given incidence relations.
A flag will be a pair (vertex, incident edge) or (vertex, incident tail). The 
set of all flags will
be denoted by $F_{\tau}$, those incident to one vertex
$v$ by $F_{\tau}(v)$.

\medskip

{\bf 1.2 Keel's presentation.}

\smallskip

Usually the cohomology ring of $\mbs$ is presented in terms of 
classes of boundary divisors as generators and quadratic 
relations as introduced by [Ke]. The additive
structure of this ring and the respective
relations can then naturally be described in terms of 
stable trees (see [KM] and [KMK]).
The boundary divisors of $\mbs$ are in one to one correspondence
with unordered 2-partitions  $\{S_1,S_2\}$ of $S$, satisfying 
$|S_1| \geq 2$ and $|S_2| \geq 2$ (stability). 
Let 
$\{ D_{\s}|\s = \{S_1,S_2\} \text{ a stable S-partition}\}$
be a set of commuting independent variables. Consider the ideal
$I_n \subset F_n$ in the graded polynomial ring
$F_S := K[D_{\{S_1,S_2\}}]$ 
generated by the following relations:
\roster 
\item "(i)" $D_{\{S_1,S_2\}} D_{\{S_1^{\prime},S_2^{\prime}\}}$,
if the number of non-empty pairwise intersection of these sets 
equal to 4.
\item "(ii)" $\forall \text{ distinct } i,j,k,l \in S: 
\sum_{ij \s kl} D_{\s} - \sum_{kj\tau il} D_{\tau}$
\endroster
where the notation of the type $ij\s kl$ is used to imply that
$\{i,j\}$ and $\{k,l\}$ are subsets of different parts of $\s$.

Set $H^*_S := F_S/I_S$.
Keels Theorem states that the map 
$$\align D_{\sigma} \longmapsto \quad &\text{dual cohomology class of the 
boundary divisor} \\
&\text{in $\overline{M}_{0n}$ corresponding to the partition $\sigma$} \endalign$$
induces the isomorphism of rings (doubling the degrees)

$$
H^*_S \overset \sim \to \longrightarrow H^* (\overline{M}_{0n},K).
\tag 1.1
$$

\medskip

{\bf 1.3 Additive structure of $H^*_n$.}

\smallskip

The additive structure of the cohomology can be nicely presented in terms of 
trees (see [KMK]). 
There (proposition 1.3) it is proved that the set of trees
with $r$ edges is in bijection with the set of good monomials of
degree $r$.
We will briefly quote some of the notions and results from that paper.
A monomial
$D_{\sigma_1} \ldots D_{\sigma_a} \in F_S$
is called good, if the family of 2-partitions
$\{\sigma_1,\ldots,\sigma_a\}$
is good, i.e.
$a(\sigma_i,\sigma_j)=3$, where for two unordered stable partitions
$\sigma =\{S_1,S_2\}$
and
$\tau =\{T_1,T_2\}$
of $S$ 
$$\align a(\sigma,\tau):= \quad &\text{the number of non-empty pairwise} \\
&\text{distinct sets among} \quad S_i \cap T_j,\, i,j =1,2. \endalign$$

\proclaim{\quad 1.3.1 Lemma (1.2 in [KMK])}
Let
$\tau$
be a stable $S$-tree with
$|E_{\tau}| \ge 1$.
For each
$e \in E_{\tau}$,
denote by
$\sigma (e)$
the 2-partition of $S$ corresponding to the one edge $S$-tree obtained by 
contracting all edges except for $e$. Then

$$mon(\tau) := \prod_{e\in E_{\tau}} D_{\sigma (e)}
\tag 1.2$$
is a good monomial.
\endproclaim

\proclaim{\quad 1.3.2 Proposition (1.3 in [KMK])}
For any
$1 \le r \le |S|-3$,
the map
$\tau \longmapsto mon(\tau)$
establishes a bijection between the set of good monomials of degree $r$ in
$F_S$
and stable $S$-trees
$\tau$
with
$|E_{\tau}| =r$
modulo $S$-isomorphism. There are no good monomials of degree greater than
$|S| -3$.
\endproclaim

{\bf 1.3.3 Additive relations.}

In [KMK] it is shown that the good monomials span the cohomology space 
and furthermore
that all linear relations between them are generated by the relative versions of (ii);
$$
\sum_{ij\tau' kl} mon(\tau') = \sum_{ik\tau'' jl} mon(\tau'')
\tag 1.3
$$
where $\{ij\tau' kl\}$ and $\{ik\tau'' lj\}$ are the preimages of the 
contraction onto a given $\tau$ contracting exactly one edge onto a 
fixed vertex $v$ seperating the flags marked by 
$i,j$ and $k,l$ resp.\ $i,k$ and $j,l$ in such a way
that they lie on different components after severing $e$, where the
markings $i,j,k,l$ refer to flags which are part of the edges of the unique 
paths from $v$ to the tails $i,j,k,l$ in $\tau$ and it required that the
paths start along different
edges.

{\bf 1.4 Trees with multiplicity.}

\sm

Since we will have to deal with monomials, which are not necessarily
good, we will extend the notion of trees to that of trees
with multiplicity.

\med

{\bf 1.4.1 Definition.}
A $S$-tree with multiplicity is a pair $(\tau,m)$ consisting of a
$S$-tree and a function $m: E_{\tau} \rightarrow \Bbb N$.

If no multiplicity function is given we will assume that it
is identically 1.

Call a monomial $D_{\sigma_1}^{m_1} \cdots D_{\sigma_k}^{m_k}$ nice if 
$a(\sigma_i,\sigma_j)= 2 \text { or }3$.

Set 
$$
mon ((\tau,m)):= \prod_{e \in E_{\tau}} D_{\s(e)}^{m(e)}.
\tag 1.4
$$

\sm

\proclaim{\quad 1.4.2 Proposition}
For any
$1 \le r \le |S|-3$,
the map:
$(\tau,m) \longmapsto mon((\tau,m))$
establishes a bijection between the set of nice monomials of degree $r$ in
$F_S$
and stable $S$-trees with multiplicity
$(\tau,m)$
with
$deg(\tau, m) := \sum_{e \in E_{\tau}} m(e) =r$.
\endproclaim

{\bf Proof.}
Immediate from 1.3.2.

\med

{\bf 1.4.3 Remark.}
Notice that unlike in the case of good monomials it can happen that a
nice monomial can represent a zero class even if the degree is less
or equal to  $|S| -3$.

\med

{\bf 1.5  Rooted trees and ordered partitions.}

\sm

{\bf 1.5.1 Remark.}
If we choose a distinguished element $s\in S$, 
we can define natural bijections between the following three sets:
\roster
\item "a)" unordered 2-partitions $\s = \{S_1,S_2\}$ of S
\item "b)" ordered 2-partitions $\s = \langle S_1, S_2\rangle$
with the condition $s \in S_2$
\item "c)" subsets $T \subseteq S \setminus \{s\}$.
\endroster
This is due to the fact that given the first component of an ordered pair of
the above type the second one is
uniquely determined. 

\pagebreak

{\bf 1.5.2 The case of $\on$.}

\sm

In particular for $S=\on$ we choose $n$ as the distinguished
element and we equivalently index the generators
of $H^*$  by subsets $S \subset \onn$ with the restriction $2 \leq |S| \leq n-2$
(note that this excludes the set $\onn$ itself).
We will denote the generator corresponding to such a set $S$:
$$D_S := D_{S,\on \setminus S} ,$$
for $S \subset \onn$.
The relations (i) and (ii) stated in this notation become:

\roster
\item "(i')" $D_S D_T$ if $S \cap T \neq \emptyset$ and the two
sets satisfy no inclusion relation.
\item "(ii')" For any four numbers $i,j,k,l$:
$$
\sum \Sb \onn \supset T \supseteq \{i,j\}\\ k,l \notin T \endSb D_T +
\sum \Sb \onn \supset T
\supseteq \{k,l\}\\ i,j \notin T \endSb D_T 
- \sum \Sb \onn \supset T^{\prime} \supseteq \{i,k\}\\ j,l \notin T \endSb 
D_T^{\prime}
- \sum \Sb \onn \supset T^{\prime} \supseteq \{j,l\}\\ i,k \notin T \endSb 
D_T^{\prime}
\tag 1.5
$$
\endroster
The expression for $D_S^2$ for a choice $i,j \in S$ and $k\notin S$ reads:
$$
D_S^2 = - \sum \Sb S \subset T \subseteq \{i,j\} \endSb D_S D_T
- \sum \Sb S \supset T \supset \onn \\ k \notin T \endSb D_S D_T.
\tag 1.6
$$
This is the formula (1.7) from [KMK] with $i,j,k,n$ playing the role
of $i,j,k,l$.

The analogs of formula (1.3) follow in the same manner.
 
\medskip

{\bf 1.5.3 Rooted trees and orientation.}

\sm

A rooted $S$-tree will be a pair $(\tau,v_{root})$ consisting 
of a $S$-tree $\tau$ and one of its vertices
$v_{root}$ called root.
An orientation of
a tree is considered to be a map 
$or:E_{\tau} \rightarrow V_{\tau}$,
with the restriction that $e$ is incident to $or(e)$. 
We will use the terminology $e$ is pointing towards $v$ to indicate
$v = or(e)$ (pointing away will be used on the same basis). 
The set $or^{-1}(v)$ will be called the incoming edges,
the remaining incident edges will be considered as outgoing.
Furthermore notice that an oriented edge $e$ of a tree
defines a subtree by cutting $e$ and selecting the tree containing $or(e)$. 
This subtree will be called the branch of $e$.

\med
 
{\bf 1.5.4 Natural orientation for a rooted tree.}

\sm

For a rooted tree $(\tau,v_{root})$ there is a natural 
orientation defined by setting
$or(e)=$ the vertex of $e$, which is furthest away from the root (i.e. $e$ is 
part of the unique path from this vertex to the root). Notice that in this 
orientation there is exactly one
incoming edge to
each vertex except for the root, which has none.
Therefore the restriction of $or$ induces an one to one correspondence of 
$V(\tau) \setminus \{v_{root}\}$ and $E(\tau)$. 
$$
\matrix \format  \r&\c&\l&\c&\r&\c&\l\\
e &\mapsto &\text {vertex to which $e$ is pointing}&
\text {\quad{\it inversely} \quad}&v &\mapsto \text {the unique incoming edge}
\endmatrix
\tag 1.7
$$

\med

{\bf 1.5.5 Orientation for an $n$-tree.}

\sm

For a given $n$-tree we will fix the root to be the vertex
with the flag numbered by $n$ emanating from it.
This defines a one to one correspondence of $n$-trees with
rooted $n$ trees.
Using this picture and remark 1.3.1 we can equivalently view a $n$-tree (with
multiplicity) as either given by the good (nice) collection of 2-partitions 
associated
to its edges or as a good (nice) collection of subsets of $\onn$ 
associated to its vertices. 
In the latter case we associate to each vertex the set $S$
of the 2-partition corresponding to the incoming edge, which 
{\it does not} contain $n$. In this way denote for given nice $\s$ and 
$S \in \s$ by  $v_S$ (resp\. $e_S$)
the vertex (resp\. edge) corresponding to $S$.

Adopting this point of view we can express quantities which are defined in the
language of Remark 1.5.1 c) in terms of  oriented $n$-trees. 
Let $\s$ be a collection
of stable subsets of $\on$, i\.e\. for each $S \in \s$ 
$S \subset \onn$ and $|S| \geq 2$.
Define for any $S \in \s$:

$$\align
\omega_{\s}(S)&= \{T| T \subset S \text{ and maximal in this respect}\}\\
depth_{\s}(S)&= |\{T| T \in \s \text{ and } T \supseteq S\}|
\tag 1.8
\endalign
$$
The definitions of (1.8) translate in the following way into tree language:

$$
\split
|S| &= |\{\text{tails marked by $i \in \onn$ on the branch of $e_S$}\}|\\
\omega_{\s}(S) &= \{\text {outgoing edges of $v_S$}\}\\
depth_{\s}(S) &= \text {the distance from $v_S$ to $v_{root}$}
\endsplit
\tag 1.9
$$
where the distance is the number of edges along the unique shortest path.

\newpage

{\bf \S 2 The intersection form}

\medskip

{\bf 2.1 Notation.}

\sm

To calculate the intersection form we need a formula for 
two monomials of complementary degree. Recall that for a tuple $(\s,m)$ 
of a nice collection of 
subsets of $\onn$ and a multiplicity function $m: \s \mapsto \Bbb N$  
we denote by
$mon(\s,m)$ the monomial $\prod_{S \in \s} D_S^{m(S)}$. 
The {\it degree} of such a
monomial is $\sum_{S \in \s} m(S)$. Furthermore let $\tau(\s,m)$ be 
the tuple $(\tau(\s), m^{\prime})$ where $\tau(\s)$ 
is the tree corresponding to the good
monomial $\prod_{S \in \s} D_S$, and 
$m^{\prime}:E_{\tau(\s)} \rightarrow \Bbb N$ is the 
the multiplicity function given by $e_S \rightarrow m(S)$.

\med

{\bf 2.2 Definition.}
A multiplicity orientation for a tree with multiplicity $(\tau,m)$ is a map 
$mult: F_{\tau} \setminus T_{\tau} \mapsto \Bbb N$ such that if $v_1$ and $v_2$ 
are the vertices of an edge $e$: 
$$mult((v_1,e)) + mult ((v_2,e))= m(e) -1. \tag 2.1$$ 
It is called good if for every $v \in V_{\tau}$ it satisfies:
$$
\sum_{f \in F_{\tau}(v)} mult(f) = |v|-3. 
\tag 2.2
$$

This is the analog of the good orientation in [KMK]. 

\proclaim {\quad 2.3 Lemma}
For a $n$-tree $(\tau,m)$ in top degree (i.e. $\sum_{e \in E_{\tau}}m(e) = n-3$)
there exists at most one good multiplicity orientation.
\endproclaim

{\bf Proof.} 
Assume that there are two good orientations $mult,mult^{\prime}$.
Consider the union of all edges on which $mult\ne mult^{\prime}.$
Each connected component of this union is a tree. Choose an end edge
$e$ of this tree and an end vertex $v$ of $e.$ At $v,$ the sum over all flags 
$f$ of $mult(f)$ and $mult^{\prime}(f)$ must be equal, but on $(v,e)$
these differ. Hence there must exist an edge $e^{\prime}\ne e$
incident to $v$ upon which $mult((v,e^{\prime}))$ and 
$mult^{\prime}((v,e^{\prime}))$
differ. But this
contradicts to the choice of $v$ and $e$.

The next lemma gives a way to decide whether this good multiplicity orientation
exists and if so to calculate it.

\proclaim {\quad 2.4 Lemma}
Assume that an $n$-tree $\tau(\s,m)$ in top degree has a good multiplicity
orientation $mult$. Let $v_S$ be the vertex corresponding to $S \in \s$ and 
$f_S$ be
the flag of the unique incoming edge then the following formula for its
multiplicity holds:
$$
mult(f_S) = |S| - 2 - \sum_{T \in \s | T \subset S} m(T).
\tag 2.3
$$
\endproclaim

{\bf Proof.}
We will use induction on the distance from the 
end vertices (i.e those vertices with only one incoming edge) in the natural 
orientation of $n$-trees
given by 1.5.5; the case for the end vertices being trivial. Now let $v_S$ be 
the vertex
corresponding to $S$. By induction we can assume that 
for all outgoing flags (2.3)
holds; i.e. for all $(v,e_T)$ with $T \in \omega_{\s}(S)$:
$$
mult((v,e_T)) = m(T)-1 - |T| +2 + \sum_{T^{\prime} \in \s | T^{\prime} 
\subset T} m(T^{\prime}). 
$$
Inserting this into the condition (2.2) we arrive at
$$
\split
mult(f_S) =& |v_S|-3 - \sum_{T \in \omega_{\s}(S)} (m(T) - 1 -mult((v,e_T)))\\
=&|S|- |\bigcup_{T \in \omega_{\s}(S)} T| + |\omega_{\s}(S)| -2 - \sum_{T \in
\omega_{\s}(S)}  (m(T) - |T| + \sum_{T^{\prime} \in \s | T^{\prime} \subset T}
m(T^{\prime}) +1) \\
=&
|S| -2 -\sum_{T \in \s | T \subset S} m(T),
\endsplit
$$
where in the last step we have used that 
$|\bigcup_{T \in \omega_{\s}(S)} T| = \sum_{T \in \omega_{\s}(S)} |T|$, 
since $\s$ is a nice
collection.

\bigskip

Consider the functional
$\int_{\overline{M}_{0,S}} : H^*(\overline{M}_{0,S}) \rightarrow K$
is given by

$$m(\tau) \longmapsto \left\{ \aligned &1, \quad \text{if} \quad \roman{deg}\ 
m(\tau) =|S|-3,\\
&0 \quad \text{otherwise.} \endaligned \right.$$
for any tree $\tau$ with $m \equiv 1$.

We put
$\langle (\tau_1,m_1) \, (\tau_2,m_2) \rangle =\int_{\overline{M}_{0S}} 
mon((\tau_1,m_1)) mon((\tau_2,m_2))$ 
and set to calculate this intersection index for the case when
$\roman{deg}\, mon((\tau_1,m_1))  +\roman{deg}\, mon((\tau_2,m_2)) =|S|-3$.
Generally, we will write
$\langle \mu \rangle$
instead of
$\int_{\overline{M}_{0S}}\mu$.

\proclaim {\quad 2.5 Theorem}
Let $mon(\s_1,m_1)$ and $mon(\s_2,m_2)$ be two monomials of complementary
degree in $H^*_n$. If there is no good multiplicity orientation of
$(\tau,m) := \tau(\s_1 \cup \s_2,m_1 +m_2)$ then
$\langle mon(\s_1,m_1) mon(\s_1,m_1) \rangle =0$. If there does exist one then:
$$
\langle mon(\s_1,m_1) mon(\s_1,m_1) \rangle 
= \prod_{v \in V_{\tau}} (-1)^{|v|-3}
\frac {(|v|-3)!} {\prod_{f \in F(v)} (mult (f))!^2}
\prod_{e\in E_{\tau}}(m(e)-1)!,
$$
where $mult$ is the unique multiplicity orientation of $(\tau,m)$ 
provided by the Lemma 2.3
whose value is given in the formula (2.3).
\endproclaim

{\bf Proof.}
Set $E:= \{e \in E_{\tau}|m(e)> 1\}$ and $\delta$ the subtree consisting 
of $E$ with
multiplicity $m\vert_{E}$
and its vertices.
Consider the canonical embedding $\varphi_{\tau}:\ 
\overline{M}_{\tau}\to \overline{M}_{0S}.$ 
$$
\langle mon(\s_1,m_1) mon(\s_1,m_1) \rangle 
=
\langle\prod_{e\in E}\varphi_{\tau}^*(D_{S(e)}^{m(e)-1})\rangle ,
\tag 2.4
$$
where the cup product in the r.h.s. is taken in
$H^*(\overline{M}_{\tau})\cong \otimes_{v\in V_{\tau}}
H^*(\overline{M}_{0,F_{\tau}(v)}).$ Applying an appropriate version
of the formulas (1.5) we can write for any
$e\in E$ with vertices $v_1,v_2$:
$$
\varphi_{\tau}^*(D_{\sigma (e)})=-\Sigma_{v_1,e} -\Sigma_{v_2,e},
\eqno{(2.5)}
$$
where 
$$
\Sigma_{v_i,e}\in H^*(\overline{M}_{0,F_{\tau}(v_i)})\otimes
\prod_{v\ne v_i}[\overline{M}_{0,F_{\tau}(v)}]
\eqno{(2.6)}
$$
and $[\overline{M}_{0,F_{\tau}(v)}]$ is the fundamental class.
Later we will choose an expression for $\Sigma_{v_i,e}$ depending
on the choice of flags denoted $i,j$ or
$k,l$ in (1.5).

\smallskip

Inserting (2.5) into (2.6), we get
$$
\langle mon(\s_1,m_1) mon(\s_1,m_1) \rangle 
=\sum_{or} \prod_{e\in E_{\tau}}(m(e)-1)!
\langle\prod_{{(v,e)\in F_{\delta}}\atop{or((v,e)) >1}} 
\frac {1}{or((v,e))!}(-\Sigma_{v,e})^{or((v,e))}\rangle,
\tag 2.7
$$
where $or$ runs over all multiplicity orientations of $\delta$. 
The summand of (2.7) corresponding to a given $or$ can be non--zero
only if for every $v\in V_{\delta}$ the sum of the degrees of factors
equals $\roman{dim}\ \overline{M}_{0,F_{\tau}(v)}=
|v|-3.$ This is what was called a good multiplicity orientation.
By Lemma 2.3 there can only exist one such orientation. 
Now assume that one good orientation $mult$ exists. We can rewrite (2.7) as
$$
\multline
\langle mon(\s_1,m_1) mon(\s_1,m_1) \rangle 
=\\
\prod_{e\in E_{\tau}}(m(e)-1)!
\prod_{{(v,e)\in F_{\delta}}\atop{mult((v,e)) >1}} 
\frac {1}{mult((v,e))!}\langle(-\Sigma_{v,e})^{mult((v,e))}\rangle .
\endmultline
\tag 2.8
$$
In view of (2.6), this expression splits into a product of
terms computed in all $H^*(\overline{M}_{0,F_{\tau}(v)}), v\in V_{\tau}$
separately. Each such term depends only on $|v|,$ and we want
to demonstrate that it equals 
$
(-1)^{|v|-3}
\frac {(|v|-3)!} {\prod_{f \in F(v)} (mult (f))!}
$. 
Put $|v|=m$, so $m\ge 3$.
Let us identify $F_{\tau}$ with $\{ 1,\dots ,m\}$ and denote by 
$D_{\rho}^{(m)}$ the class of a
boundary divisor
in $H^*(\overline{M}_{0,m})$ corresponding to a stable
partition $\rho$ of $\{1,\dots ,m\}$ and set $d_i := mult((v,e_i))$, 
where $e_i$ is the edge
belonging to the flag $i\in \{1,\dots ,m\}$. The contribution of $v$ 
in (2.8) becomes
$$
\prod_{i=1}^m \langle (-\Sigma_i^{(m)})^{d_i} \rangle 
:= g(d_1,\dots,d_m),
\tag 2.9
$$
where $-\Sigma_i^{(m)}$ is the element of (2.6) and the superscript $(m)$ 
is again included
to keep track of the spaces involved.
We will prove the following properties of the function $g(d_1,\dots,d_m)$ 
identifying it as
$(-1)^{m-3}\frac{(m-3)!}{d_1!\dots d_m!}$.
\roster
\item "a)" $g(0,0,0) = 1$.
\item "b)" $g(d_1,\dots,d_m)$ is symmetric in the $d_i$.
\item "c)" If $d_m = 0$ then 

$g(d_1,\dots,d_m) = -\sum_{i: d_i >1} g(d_1,\dots,d_i-1,\dots,d_m)$.
\endroster

\med

{\bf 2.5.1 Remarks.}
Notice that up to the minus sign in c) these are exactly the conditions 
satisfied by the
numbers $\langle \tau_{\alpha_1} \dots \tau_{\alpha_m} \rangle$
in genus zero [K].
Furthermore we can always choose the flags in such a way that the flags 
$1,\dots ,k$ ($k \leq m-3$) belong to the edges $e$ with 
$mult(f(v,e)) >1.$

\roster
\item "ad a)" We have by definition 
$\langle[ \overline{M}_{0,3}]\rangle =1$.
\item "ad b)" The symmetricity results from the fact that the integral 
in question does not
depend on a renumbering of the flags.
\item "ad c)"
First we can use relation (2.5) for any $k,l$ to write
$$
- \Sigma_i^{(m)} = \sum_{\rho:\ i\rho\{k,l\}}-D^{(m)}_{\rho}
\tag 2.10
$$
We will calculate (2.9) inductively. Consider the projection map
(forgetting the $(m)$--th point) $p:\ \overline{M}_{0,m}\to
\overline{M}_{0,m-1}$ and the $i$--th section map
$x_i:\ \overline{M}_{0,m-1}\to\overline{M}_{0,m}$
obtained via the identification of $\overline{M}_{0,m+3}$ with the
universal curve. We have $p\circ x_i=\roman{id},$ and 
$x_i$ identifies $\overline{M}_{0,m-1}$ with $D^{(m)}_{\sigma_i}$
where
$$
\sigma_i = \{\{m,i\}\{1,\dots,\widehat{i},\dots,m-1\}\}; 
$$
so if we choose some $k,l \neq m$:
$$
\sum_{\rho:\ i\rho\{k,l\}}-D^{(m)}_{\rho}=
-p^*\left(\sum_{\rho^{\prime}:\ i\rho^{\prime}\{k,l\}}
D^{(m-1)}_{\rho^{\prime}}\right)
-x_{i*}([\overline{M}_{0,m-1}]).
\eqno{(2.11)}
$$
We will now replace one of the $\Sigma_i$ for each $i$ with $d_i >1$ 
using (2.10)
with some arbitrary $k,l \neq m$. Then (2.9) reads
$$
\prod_{i=1}^m \langle
\left( -p^*(\sum_{\rho^{\prime}:\ i\rho^{\prime}\{k,l\}}
D^{(m-1)}_{\rho^{\prime}})
-x_{i*}([\overline{M}_{0,m-1}]) \right)
(-\Sigma_i^{(m)})^{d_i-1}\rangle
\eqno{(2.12)}
$$
where $\rho^{\prime}$ runs over stable partitions of
$\{1,\dots ,m-1\}.$ We represent the resulting expression as a sum of products
consisting of several $p^*$--terms and several $x_{i*}$--terms each.
If such a product contains $\ge 2\ x_{i*}$--terms, it
vanishes because the structure sections pairwise do not intersect.
We obtain
$$
\split
\sum_{i:d_i >1} \langle \prod_{j \neq i:d_j >1} \left(-p^*(\sum_{\rho^{\prime}:\
j\rho^{\prime}\{k,l\}}
D^{(m-1)}_{\rho^{\prime}}) (-\Sigma_j^{(m)})^{d_j-1} \right)
(-x_{i*}([\overline{M}_{0,m-1}]))\rangle\\
+ \langle \prod_{i:d_i>1} p^*(-\sum_{\rho^{\prime}:\ i\rho^{\prime}\{k,l\}}
D^{(m-1)}_{\rho^{\prime}}) (-\Sigma_j^{(m)})^{d_j-1} \rangle .
\endsplit
\tag 2.13
$$
If $d_i-1 >0$ then the summand containing an $x_{i*}$--term will vanish. 
To see this again
replace one of the $\Sigma_i$ using
(2.10) but with $k=m$ and some $l$. In case $d_i-1 =0$ we can write 
the respective term in the sum in (2.13) as
$$
\langle (p^*(-\sum_{\rho^{\prime}:\ j\rho^{\prime}\{k,l\}}
D^{(m-1)}_{\rho^{\prime}})^{d_j-1})(-x_{i*}([\overline{M}_{0,m-1}]))
\rangle
$$
by replacing the $\Sigma_j$ according to (2.11) and again using the
fact that the structure sections do not pairwise intersect. 
Using induction on the last summand in (2.13) we arrive at the situation, where 
all $\Sigma_i^{(m)}$'s have been replaced. And the product only contains 
$p^*(\Sigma_i^{(m-1)})$-term but this term vanishes because
$\roman{dim}\ \overline{M}_{0,m-1}=m-2$.
Finally, we are left with for one summand for each $i:d_i >1$ containing 
only one
$x_{i*}$--term and $p^*$--terms. Using the projection formula
$$
\langle p^*(X) x_{i*}([ \overline{M}_{0,m+2}])\rangle=\langle X\rangle
$$
one sees that each such term equals $-g(d_1,\dots d_i-1,\dots,d_{m-1})$. 
And the result follows.

\endroster

\newpage

{\bf \S 3  A boundary divisorial basis and its tree representation}

\medskip

The work presented in this section is inspired by the 
presentation of a basis of the cohomology ring of $\mbn$ given in terms of 
hyperplane
sections in [Yu]; especially the notions of the $^*$-operation and the 
order have been adapted to the present context. 

\med

{\bf 3.0 Preliminaries.}

\sm

In order to state the basis we make use of certain classes
$$
D_S x_S^k := \pi_{f_S*}(D_S^{k+1} D_{S\amalg {f_S}}), \quad k \geq 0
\tag 3.1
$$
where $\pi_{f_S*}:\overline{M}_{0, \on\amalg f_S} \rightarrow  \mbn$
is the forgetful map forgetting the point $f_S$. 

Another way to
present these classes is given by the following observation. 
Consider the following decomposition of $D_S^2$ using (1.6):
$$
D_S^2 = D_S (\sum_{\{i,j\} \subset T \subset S} D_T + 
\sum\Sb \onn \supset T' \supset S \\ k \notin T' \endSb D_T') =: D_S(x_S + y_S)
\tag 3.2
$$
for any choice of $i,j \in S, k,l \notin S$. 
With the notation (3.2) we can write $D_S^{k+1}$ in the same spirit as:
$$
D_S^{k+1} = D_S(\sum_{i=0}^k \binom {k}{i} x_S^i y_S^{k-i}).
\tag 3.3
$$
In the context of the proof of theorem 2.5 each summand of (3.3) corresponds to a choice of multiplicity orientation.
In particular the term with $x_S^i$ corresponds to the one which satisfies
$mult (f_S)= i, mult(f_{S^c})= k-i$, for the flags $f_S$ and $f_{S^c}$ of
$e_S$ so that we
can identify (3.1) with the summand corresponding to 
$mult(f_S)=k, mult (f_{S^c})=0$.

\med

{\bf 3.0.1 A tree representation.}

\sm

A tree representation for a class (3.1) is given by a choice
an ordered $k+1$ element subset $\langle f_1, \dots, f_{k+1}\rangle$ of $S$ as 
the sum over all assignments
of the flags of \linebreak
$S \setminus \{f_1, \cdots, f_{k+1}\}$ to the vertices
of the linear tree determined by \linebreak
$D_{\{f_1,f_2\}} D_{\{f_1,f_2,f_3\}} \dots D_{\{f_1, \dots, f_{k+1}\}}$

$$
D_S x_S^k = D_S \sum \Sb \langle S_1, \dots , S_k \rangle \\ S_1 \amalg \dots
\amalg S_k = S \setminus \{f_1, \dots, f_{k+1}\} \endSb 
D_{\{f_1,f_2\} \amalg S_1}
D_{\{f_1,f_2,f_3\}\amalg S_2} \dots D_{\{f_1, \dots, f_{k+1}\}\amalg S_k}
\tag 3.4
$$
 
or more generally let  $\tau$ given by $D_{T_1} \cdots D_{T_k}$ be any tree
with $|v_{T_i}| = 3$ for $i= 1,\dots,k$ and $T_1 \amalg \dots \amalg T_k =
\{f_1, \dots , f_k\}$ then 

$$
D_S x_S^k = D_S \sum \Sb \langle S_1, \dots , S_k \rangle \\ S_1 \amalg \dots
\amalg S_k = S \setminus \{f_1, \dots, f_{k+1}\} \endSb D_{T_1 \amalg S_1}
\dots D_{T_k \amalg S_k}.
\tag 3.5
$$

Both (3.4) and (3.5) follow from (1.5) with the appropriate choices for
the flags.

\medskip

{\bf 3.1 The basis.}

\sm

Consider a class of the following type
$$
\mu = \pi_{n*} (D_{S_1} x_{S_1}^{m(S_1)} \cdots D_{S_k} x_{S_k}^{m(S_k)} 
D_{\onn} x_{\onn}^{m(\onn)}), \quad m(S) \geq 0
\tag 3.6
$$
To this class we associate the underlying n-tree $\tau(\mu)$ determined by
$D_{S_1} \dots D_{S_k} D_{\onn}$.  The powers $m(S)$ then determine a unique 
multiplicity orientation in the sense of 3.0 given by 
$mult (f_S) := m(S), mult (f_{S^c}) = 0$, where $f_S$ and $f_{S^c}$ are the 
flags corresponding
to the edge $e_S$ in $\tau(\mu)$.

Using the equations of the type (3.4) we can associate to each monomial $\mu$
a sum of good monomials, which we will call $tree (\mu)$.

Consider the following set

$$
\multline
{\Cal B}_n := \{ \pi_{n+1*} (D_{S_1} x_{S_1}^{m(S_1)} \cdots D_{S_k} 
x_{S_k}^{m(S_k)} D_{\onn} x_{\onn}^{m(\onn)}) \,
\vert \, 0 \leq  m(S) \leq v_S -4 \text { and } \\
0 \leq m(\onn) \leq v_{\onn}-3\}
\endmultline
\tag 3.7
$$

\proclaim {3.1.1 Proposition}
The set ${\Cal B}_n$ is a basis for $A^*(\mbn)$.
\endproclaim

{\bf Proof.} 
By Lemma 3.1.2 and 3.1.5.

\proclaim {3.1.2 Lemma} The set ${\Cal B}_n$ spans $A^*(\mbn)$.
\endproclaim

{\bf Proof.} From [Ke] and [KMK] we know that the good monomials span,
so it will be sufficient to show that any such monomial is in the span of 
${\Cal B}_n$. Now let $\tau(\mu)$ be the tree corresponding to such a 
good monomial $\mu$.
If for all $v \in V_{\tau}\, |v| \geq 4$, then the monomial is already in 
${\Cal B}_n$.
If not let $\tau_3$ be a maximal subtree of $\tau$ whose vertices
except for the root (induced by the natural orientation) all have valency 
three; call such a tree a  3-subtree and the number of its edges
its length.  Furthermore let $R$ be the set associated with the root. 
Let $F_3(\tau_3)$ the set of tails of $\tau_3$ without the ones coming from 
the root.
The formula (3.5) for the tree representation of $D_R x_R^{l}$ with the choice
of $F_3(\tau_3)$ as the fixed set and $\tau_3$ as a 3-subtree expresses 
$\tau$ in terms of trees with less maximal 3-subtrees of maximal length,
whose vertices either comply with the conditions of ${\Cal B}_n$ or 
are part of a unique maximal whose root $v_R$ has multiplicity $0$, 
i.e. $x_R$ does not divide the monomial corresponding to the tree. 
Notice that if the root $v_R$ of any 3-subtree is three valent then $R= \onn$.
We can now proceed by induction of the number of such maximal 3-subtrees 
with the maximal number of edges $l$.

\med

{\bf 3.1.2 The *-operation.}

\sm

We define the following involution on ${\Cal B}_n$:

$$
\multline 
\pi_{n+1*} (D_{S_1} x_{S_1}^{m(S_1)} \cdots D_{S_k} x_{S_k}^{m(S_k)} D_{\onn}
x_{\onn}^{m(\onn)}) 
\; {\overset * \to \longrightarrow}\\
\pi_{n+1*} ((-1)^{|v_{S_1}|-3} D_{S_1} x_{S_1}^{|v_{S_1}| -4 - m(S_1)} \cdots
(-1)^{|v_{S_k}|-3} D_{S_k} x_{S_k}^{|v_{S_k}|-4-m(S_k)}\\
(-1)^{|v_{\onn}|-3}D_{\onn}x_{\onn}^{|v_{\onn}| -3 - m(\onn)}).
\endmultline
\tag 3.8
$$

This operation preserves the underlying tree $\tau(\mu)$ but changes the 
multiplicities in such a way that $\mu$ and $\mu^*$ have complimentary 
dimensions. More precisely consider $\mu$ as the push forward of the class 
$\bigotimes_{v_S \in V_{\tau(\mu)}} x_S^{m(S)} \in 
H^*(\overline {\Cal M}_{\tau(\mu)})$ to 
$H^*(\mbn)$, then locally at each vertex we have a class of degree $m(S)$.
This class is replaced under the *-operation by a ``dual'' class of 
complimentary 
degree $\dim (\ol{M}_{0,F_{\tau}} (v_S)) -m(S)$, which is provided as a summand
of $\vp_{D_S}^*(D_S x_S^{|v_S|-4-m(S)})$.

\proclaim {3.1.3 Lemma} 

For two elements $\mu,\nu$ of ${\Cal B}_n$ the integral $\int_{\mbn} \mu \nu^*$ 
does not vanish iff $\tau(\mu\nu^*)$ is nonzero and if there is one good 
multiplicity orientation among the multiplicity orientations
satisfying  $(f_S) = m^{\mu}(S) + m^{\nu^*}(S) +1, mult (f_{S^c}) = 0$ or  
$(f_S) = m^{\mu}(S) + m^{\nu^*}(S), mult (f_{S^c}) = 1$, where $f_S,f_{S^c}$ 
are the flags of the edge $e_S$. If such an orientation exists it is unique 
and 
$$
\int_{\mbn} \mu \nu^*
 = \prod_{v \in V_{\tau (\nu)}} (-1)^{|v|-3} \prod_{v\in V_{\tau(\mu\nu^*)}} 
(-1)^{|v|-3} \frac {(|v|-3)!}
{\prod_{f\in F_{\tau(\mu\nu^*)}(v)} (mult (f))!}.
\tag 3.9
$$

\endproclaim

{\bf Proof.}
The formula (3.9) and the conditions for $\mu$ and $\nu$ as well as
the ones for the considered multiplicity orientations follows from theorem 2.5
by considering the summands of 
$$
\pi_{n+1*}
(D_{S_1}^{\epsilon (S_1)+m(S_1)} \cdots D_{S_l}^{\epsilon (S_l)+m(S_l)} 
D_{\onn}^{m(\onn)}).
$$
corresponding via 3.0 to the given monomial
$$\mu \nu^*= \pi_{n+1*}
(D_{S_1}^{\epsilon (S_1)} x_{S_1}^{m(S_1)} \cdots D_{S_l}^{\epsilon (S_l)} 
x_{S_l}^{m(S_l)} D_{\onn} x_{\onn}^{m(\onn)})$$  
with $\epsilon (S) \in \{1,2\}$.

Notice that in the formula (3.9) the  binomial coefficients
$\binom {m(e_S)-1}{mult(f_S)}$ which appear in theorem 2.5 are absent.
This is due to the fact that these factors stemming from the expansion of
$D_S^{m(e_S)}$ as in (3.3) are stripped off in the definition of the classes
$D_S x_S^k$. 

\med

{\bf 3.1.4 An order.}

\sm

Given two monomials $\mu, \mu'$ of type (3.6) of the same degree we call 
$\mu \prec \mu'$ if for the maximal integer
$k$ such that all sets of the depth $d$ vertices for $1\leq d \leq k$ 
coincide and
$m(S)=m'(S)$ for all sets of the depth $d'$ vertices, for $1\leq d' < k$ 
one of the following conditions
holds
\roster
\item "(a)" $m(S)\leq m'(S)$ for all $S$ of depth $k$  and the
inequality is strict for at least one $S$ or
\item "(b)"  $m(S) = m'(S)$  and  $|v_S| \leq |v'_S|$ for all $S$ of depth $k$
and there is at least one $S$ where the inequality is strict.
\endroster 

It is easy to check that this defines a half order on ${\Cal B}_n$.

The *-operation connects with the half order $\prec$ in the
following way:

\proclaim {3.1.5 Lemma} If $\mu,\nu \in {\Cal B}_n$ are two distinct
basis elements ($\mu \neq \nu$) and $\mu \nu^* \neq 0$ then $\mu \prec \nu$.
\endproclaim

{\bf Proof.}
We will use superscripts $\mu,\nu$ to refer to the quantities
concerning the monomials $\mu,\nu$ and take quantities without any superscript 
to refer to $\mu\nu^*$. So the notation $|v_S^{\nu}|$ is used 
for the valency of the vertex $v_S$ in
the tree $\tau(\nu)$ and  $|v_S|$  without any superscript is taken to be 
the valency of the vertex $v_S$ in the tree $\tau(\mu \nu^*)$. If 
$\mu \nu^* \neq 0$ then the underlying tree of $\mu \nu^*$ carries a unique
good multiplicity orientation by theorem 2.5. Furthermore the underlying 
trees of $\mu$ and $\nu$ coincide up to depth $k$; this is the 
first condition for $k$. From this together with Lemma 3.1.3 it 
follows that the good multiplicity orientation up to depth $k-1$ 
is given by
$mult(f_S) = |v_S| -3$.  Now at depth $k$ we
must have $mult(f_S) \leq |v_S|-3$ and because the multiplicity orientation is
fixed for all lower depths as specified we also have
$mult(f_S) = m^{\mu}(S) + m^{\nu^*}(S) + \delta_{S,\onn} = 
m(S) + |v_S^{\nu}|-3- m^{\nu}(S)$. Combining these two relations 
we find the condition:
$$
m^{\mu}(S) - m^{\nu}(S) \leq |v_S| - |v_S^{\nu}|.
\tag 3.10
$$
Furthermore we have the inequalities $|v_S| \leq |v_S^{\nu}|, |v_S| 
\leq |v_S^{\mu}|$, since $\tau(\mu)$ and $\tau(\nu^*)=\tau(\nu)$ result 
from $\tau(\mu \nu^*)$ via contractions of edges which only increase the 
number of flags at the remaining vertex. So the left hand side of
(3.10) is less or equal to zero:
$$
m^{\mu}(S) - m^{\nu}(S) \leq 0.
\tag 3.11
$$ 
Thus if the inequality is strict for some $S$ we arrive at condition (a), if 
however $m^{\mu}(S)=m^{\nu}(S)$ for all $S$ of depth $k$ the 
following inequality must also hold:
$$
0 \leq |v^{\mu}(S)| - |v^{\nu}(S)|.
\tag 3.12
$$
Equality for all $S$ in (3.12) however would contradict the choice of $k$ since
if $m^{\mu} (S)=m^{\nu} (S)$ and $|v^{\mu}(S)| = |v^{\nu}(S)|$ we have 
$|v^{\nu}(S)| = |v(S)|=|v^{\mu}(S)|$ from the above inequalities, so that 
there are no contractions from $\tau(\mu\nu^*)$ to $\tau(\mu)$ and $\tau(\nu)$
up to depth $k+1$ and the sets of depth $k+1$ corresponding
to the outgoing edges of $v^{\mu}$ and  $v^{\nu}(S)$ must also coincide.

\proclaim {3.1.6 Lemma} 

Consider the matrix $T=(t_{\mu,\nu})_{\mu,\nu \in {\Cal B}_n}$  given by
$$
t_{\mu,\nu} := \int_{\mbn} \mu \nu^*.
$$ 
This matrix is unipotent and the entry $t_{\mu,\nu}$ is determined by 
Lemma 3.1.3.

In particular, the set $\Cal B$ is linear independent.
\endproclaim

{\bf Proof.}
For the diagonal entries $\int \mu \mu^*$ $mult(f_S)=|v_S|-3$ is a good 
multiplicity orientation so that (3.9) renders $t_{\mu\mu^*}=1$. 
Furthermore by considering any extension of the half order to a total order
the unipotency is proved by Lemma 3.1.5. 

\pagebreak

{\bf 3.2 The intersection form and its inverse for the basis ${\Cal B}_n$.}

\sm

With the help of the matrix $T$ introduced in 3.1.6 we can write the matrix 
$M$ for the intersection form in the basis ${\Cal B}_n$ as  $M = TP$, where 
the matrix  $P$ is the matrix representation of the $^*$-operation given by 
the signed permutation matrix 
$$
P_{\mu, \nu} = (-1)^{n-3-|E_{\tau(\mu)}|} \delta_{\mu, \mu^*}
\tag 3.13
$$. 

\proclaim {Theorem 3.2.1} The Gram-matrix $(m_{\mu\nu})$ for the basis 
${\Cal B}_n$ is given by
$$
m_{\mu\nu} = (-1)^{n-3-|E_{\tau(\nu)}|} t_{\mu \nu^*}
\tag 3.14
$$

and its inverse matrix $(m^{\mu \nu})$ is given by the formula:

$$
m^{\mu \nu} = (-1)^{n-3-|E_{\tau(\mu)}|} (\delta_{\mu^* \nu} + \sum_{k\geq 0} 
\sum_{\mu^* \prec \tau_1 \dots \prec \tau_k \prec\nu} t_{\mu^* \tau_1} t_{\tau_1 
\tau_2} \cdots t_{\tau_{k-1} \tau_{k}}
t_{\tau_{k} \nu})
\tag 3.15
$$
where the values for the $t_{\sigma,\sigma'}$ are given by (3.9) 
and the sum over 
$k$ is finite. 
\endproclaim
 
{\bf Proof.}

Formula (3.14) follows from the above decomposition $M=TP$. 
To prove formula (3.15)
set $N:= id-T$. According to Lemma 3.1.6 $N$ is nilpotent and the 
inverse to the intersection form can now be written as 
$$
M^{-1}= PT^{-1}= P(id + N + N^2 + \dots)
\tag 3.16
$$ 
where the sum in (3.16) is finite.

\newpage

{\bf \S 4 Applications to Frobenius manifolds and quantum cohomology}

\med

{\bf 4.1 Particular cases.}

\sm

Writing down the results of \S 2 and \S 3 we obtain the following intersection matrices $M_n$ for small values of $n$:
\roster
\item "$n=3$" $M_3 = (1)$.
\item "$n=4$" For the basis $\pi_{5*}(D_{1,2,3}), \pi_{5*}(D_{1,2,3}x_{1,2,3})$ we obtain
$$M_4 = \left ( \matrix
0& 1\\ 1&0 \endmatrix \right )$$ 

\item "$n=5$" For the basis $\pi_{6*}(D_{1,2,3,4}), D_{1,2,3}, 
D_{1,2,4},D_{1,3,4}, D_{2,3,4}, \pi_{6*}(D_{1,2,3,4}x_{1,2,3,4})$, \linebreak 
$\pi_{6*}(D_{1,2,3,4}x_{1,2,3,4}^2)$ the intersection matrix is:
$$M_5 = \left ( \matrix 0&0&0&0&0&0&1\\  0&-1&0&0&0&0&0\\ 0&0&-1&0&0&0&0\\
0&0&0&-1&0&0&0\\ 0&0&0&0&-1&0&0\\ 0&0&0&0&0&1&0\\ 1&0&0&0&0&0&0
\endmatrix 
\right ) $$

\item "$n=6$" In this case the intersection matrix also has only nonzero entries
for the integrals of dual classes under the $*$-operation: 
$\int_{\mbn} \mu\mu^*$ whose values are $(-1)^{3-|E_{\tau(\mu)}|}$.

\item "$n\geq 7$"  For the higher values of $n$ the structure of the matrix $T$ 
is not diagonal since also entries other than those coming from the product of 
*-dual classes  can be nonzero e.g. 
$\langle D_{i,j,k,l}x_{i,j,k,l} D_{i,j,k,l}x_{i,j,k,l} \rangle$ in 
${\overline M}_{0,7}$. Thus the *-operation fails to give the Poincar\'e
duality for these spaces.

However on the subspace ${A^1(\mbn) \oplus A^{n-4}(\mbn)}$ the  *-operation 
does provide the Poincar\'e duality as can be deduced from Lemma 3.1.3. 
On this subspace the matrix $T$ is just the
identity matrix, so that the restriction to this subspace of $M_n$ is given by 
$P$. In the case of small $n <7$ this subspace is already the whole space, so 
that the matrices in the previous cases are just given by $P$.

\endroster

\med

{\bf 4.2 Tensor product of higher order operations of formal Frobenius 
manifolds.}

\med

{\bf 2.4.1 Frobenius manifolds}

\sm

A formal Frobenius manifolds is a triple $(H,g,$ additional structure), where 
$H$ is a (super) vector space over a field $K$ of characteristic zero, $g$ is
a non-degenerate even scalar product on $H$ and the additional structure is one 
of the following [D, KM, KMK]: 
\roster
\item "(i)" a Cohomological Field Theory (CohFT) ($I_n$),

\item "(ii)" a potential $\Phi$ for a 
set of  Abstract Correlation Functions (ACF) ($Y_n$) satisfying the 
WDVV-equations or 

\item "(iii)" a structure of cyclic $C_{\infty}$-algebra on $H$ ($\circ_n$). 
\endroster

The moduli space of rank one CohFT and 
the respective structure of tensor product is presented in [KMK] and [KMZ].

As a brief reminder we recall that a CohFT on $(H,g)$ is  given by
a series of ${\Bbb S}_n$-equivariant maps:
$$
I_n: H^{\otimes n} \rightarrow H^*(\mbn,K), \quad n\geq 3
$$
which satisfy the relations:
$$
\varphi^*_{\sigma} (I_n (\gamma_1 \otimes \ldots \otimes \gamma_n)) =\epsilon 
(\sigma) (I_{n_1+1} \otimes I_{n_2+1}) (\bigotimes_{j\in S_1} \gamma_j \otimes 
\Delta \otimes (\bigotimes_{k \in S_2} \gamma_k))
\tag 4.1
$$
where $\vp_{\sigma}$ for $\sigma = S_1 \amalg S_2$ is the inclusion map of the 
divisor $D_{\sigma}$, 
$\vp_{\sigma}: \ol{M}_{0,|S_1|+1} \times \ol{M}_{0,|S_2|+1} \rightarrow \mbn$,
$\Delta =\Sigma \Delta_a \otimes \Delta_b g^{ab}$
is the Casimir element, and
$\epsilon (\sigma)$
is the sign of the permutation induced on the odd arguments
$\gamma_1, \ldots, \gamma_n$.

\med

{\bf 4.2.2 Equivalences of the different structures}

Given a CohFT the associated system of ACF'a is defined as follows:

$$
Y_n (\gamma_1 \otimes \dots \otimes \gamma_n) = \int_{\mbn}I_n (\gamma_1 \otimes
\dots \otimes \gamma_n).
\tag 4.2
$$
The potential for a system of ACF's is given (after a choice of a 
basis $\{\Delta_a\}$ and dual coordinates $x^a$ of $H)$ as a formal power series
depending on a point $\gamma = \sum x^a D_a$ by:
$$
\Phi(\gamma) = \sum_{n\geq 3} \frac{1}{n!} Y_n((x^a \Delta_a)^{\otimes n}).
\tag 4.3
$$
The conditions (4.1) on the $I_n$ are equivalent to the WDVV or
associativity equations [KM]:
$$
\sum_{ef} \partial_a \partial_b \partial _c \Phi \cdot g^{ef} \partial_f 
\partial_c \partial_d \Phi 
=(-1)^{a(b+c)} \sum_{ef} \partial_b \partial_c \partial_e \Phi \cdot g^{ef} 
\partial_f \partial_a\partial_d \Phi. 
\tag 4.4
$$
Here we use the simplified notation
$(-1)^{a(b+c)}$
for
$(-1)^{\tilde x_a(\tilde x_b +\tilde x_c)}$
where
$\tilde x$
is the
${\Bbb Z}_2$--degree of $x$. 

The reverse direction of (4.2), i.e.\ the reconstruction of a CohFT from its
ACF's is contained in 
the second reconstruction theorem of [KM]. 
In this context the $I_n$ can be recovered by extending the $Y_n$ to a set
of operadic ACF, i.e. a set $\{Y(\tau)| \tau$ is a $n$-tree$\}$ satisfying  
$$
Y(\tau)(\gamma_1 \dots \gamma_n) = (\bigotimes_{v\in V(\tau)} Y_{|v|})
(\gamma_1 \otimes \dots \otimes \gamma_n \otimes \Delta^{\otimes E_{\tau}})
\tag 4.5
$$
where the arguments $Y_{|v|}$ are labeled by the flags of $v$ (for the precise
formalism of operadic ACF see [KM]).
the $I_n$ themselves can be calculated via their Poincar\'e duals 
with the help of the formula:
$$
Y(\tau)(\gamma_1\otimes \ldots \otimes \gamma_n)
=\int_{\ol{M}_{\tau}} \vp^*(I_n (\gamma_1\otimes \ldots \otimes \gamma_n)).
\tag 4.6
$$

The explicit calculation of the maps $I_n$ given potential
$\Phi$ or a set of $Y_n$ thus depends
on the knowledge of the Poincar\'e duality as noted in [KMK] and
is made possible by the results of \S 2 and \S 3.

The higher order multiplications are derived from the ACF's in the
following manner:
$$
\circ_n:= H^{\otimes n} @> Y_{n+1} >> \check{H} @> g >> H
\tag 4.7
$$

In the operadic setting given a set of higher order multiplication 
there is a natural operation associated to each $n$-tree $\tau$ (see [GK])
which we will call $\circ(\tau)$. Such an operation corresponds to a cyclic 
word with parenthesis roughly as follows. Denote the multiplication $\circ_n$
by the word $(x_1, \dots, x_n)$ and think of it as a one vertex tree with $n$
incoming flags and one outgoing flag. Composing two higher multiplications 
corresponds to grafting two such trees in such a way that the outgoing flag of
one tree is fused together with one of the incoming flags
of the other tree to form an edge, e.g.\ the flag $i$ for 
$(x_1,\dots,x_{i-1},(x_i, \dots , x_{i+k}), x_{i+k+1}, \dots, x_n)$. 
Continuing in this way we obtain a tree from any such word
and vice versa we can 
associate to any $n$ tree with the orientation of 1.5.5 a $(n$--$1)$-ary 
operation
of composed higher multiplications.

\med

{\bf 4.2.3 Tensor product for Frobenius manifolds}

\sm

In the language of CohFT the tensor product of two formal Frobenius
$(H',g',\{I'_n\})$ and $(H'',g'',\{I''_n\})$
is given by the tensor product CohFT on $H' \otimes H''$  which is naturally 
defined via the cup product in $H^*(\mbn,K)$:

$$ (I'_n\otimes I''_n)(\gamma'_1 \otimes \gamma''_1 \otimes \dots \otimes 
\gamma'_n
\otimes \gamma''_n) := \epsilon (\gamma',\gamma'') I'_n(\gamma'_1 \otimes \dots 
\otimes \gamma'_n) \wedge
I''_n(\gamma''_1 \otimes \dots \otimes \gamma''_n)
\tag 4.8$$
where $\epsilon (\gamma',\gamma'')$ is the superalgebra sign.

Using (4.2 - 4.7) one can transfer this definition of the tensor product onto
any of the other structures $(Y_n, Y(\tau),\Phi, \circ_n, \circ(\tau))$. 

In particular using
$Y'_n$ and $Y''_n$ we obtain:
$$
\multline
(Y'_n \otimes Y''_n)(\gamma'_1 \otimes \gamma''_1 \otimes \dots \otimes 
\gamma'_n
\otimes \gamma''_n) = \\
\int_{\mbn}I'_n(\gamma'_1 \otimes \dots \otimes \gamma'_n) \wedge
I''_n(\gamma''_1 \otimes \dots \otimes \gamma''_n)\\
\endmultline
\tag 4.9
$$

In order to calculate the integrals on the right hand side of (4.5) we will
use the basis, the calculation of its Gram matrix and its inverse obtained 
in previous paragraph. We also utilize the operadic 
correlation functions corresponding to $Y'_n,Y''_n$ (see [KM]) and use the
notation $Y(\mu)$ for $Y(tree(\mu))$ for a $\mu$ in ${\Cal B}_n$. 
Now, let ${\Cal B}_n$ be the basis of $H^*(\mbn)$ given in 3.1 and
$\check{\mu} = \sum_{\mu\nu} m^{\mu \nu} \nu$ the dual basis. Combining
the results of \S 3 with the formula (4.6) we obtain the following:

\pagebreak

\proclaim {4.2.4 Corollary} The tensor product of two CohFT $(H',g',Y')$ and
$(H'',g'',Y'')$ is given by:

$$\multline
(Y'_n \otimes Y''_n)(\gamma'_1 \otimes \gamma''_1 \otimes \dots \otimes 
\gamma'_n
\otimes \gamma''_n) = \\
\sum_{\mu, \nu \in {\Cal B}_n} 
Y'(\check {\mu}) (\gamma'_1 \otimes \dots \otimes \gamma'_n) m_{\mu \nu}
Y''(\check{\nu}) (\gamma''_1 \otimes \dots \otimes \gamma''_n).
\endmultline
\tag 4.10
$$ 
\endproclaim

\proclaim {4.2.5 Corollary} The tensor product of two Frobenius manifolds in
terms of the higher order multiplications
is given by
$$
\multline
\circ'_n \otimes \circ''_n (\gamma'_1 \otimes \gamma''_1 \otimes \dots \otimes 
\gamma'_n
\otimes \gamma''_n) = \\
\sum_{\mu, \nu \in {\Cal B}_n} 
\circ'(\check {\mu}) (\gamma'_1 \otimes \dots \otimes \gamma'_n) m_{\mu \nu}
\circ''(\check{\nu}) (\gamma''_1 \otimes \dots \otimes \gamma''_n).
\endmultline
\tag 4.11
$$
\endproclaim

{\bf 4.3 The K\"unneth formula in quantum cohomology}

\med

{\bf 4.3.1 Quantum cohomology}

\sm

The quantum cohomology of a projective manifold $V$ will be regarded as
a formal deformation of its cohomology ring with the coordinates of the space 
$H^*(V)$ being the
parameters. The structure constants are given by a formal series $\Phi^V$,
which is defined in terms of Gromov-Witten invariants [KM]. One can regard the 
quantum cohomology as the a structure of Frobenius manifold on $(H^*(V)$, 
Poincar\'e pairing) with the GW-invariants playing the role of the $I_n$ 
and the potential $\Phi^V$
being the potential of (4.2). 
The quantum cohomology of a product $V \times W$ regarded
as a Frobenius manifold is just the tensor product of the Frobenius
manifolds $H^*(V)\otimes H^*(W), \text {Poincar\'e pairing}, \Phi^{V\times W}$,
as can be shown using [B]\footnote{K.\ Behrend private communication}.
Putting together (4.4) and the corollary  4.2.1 we obtain the explicit K\"unneth
formula:

\proclaim {Corollary 4.3.2}
The potential $\Phi^{V\times W}$ of the quantum cohomology of
$V \times W$ is given by the formula:
$$
\Phi^{V\times W}(\gamma' \otimes \gamma'') = \sum_{n\geq 3} 
\sum_{\mu, \nu \in {\Cal B}_n} 
Y'(\check {\mu}) (\gamma^{\prime \otimes n} ) m_{\mu \nu}
Y''(\check{\nu}) (\gamma^{\prime \prime \otimes n}).
\tag 4.12
$$
\endproclaim

\med

{\bf 4.4 Examples.}

\med

{\bf 4.4.1 Higher order correlation functions.}

\sm

Using the calculations of 4.1 we obtain the following formulas for the tensor
product of the first higher order correlation functions of $(H',g',Y'_n)$ and
$(H'',g'',Y''_n)$. To write down the
formulas let $\sum_{a'b'} \Delta_{a'} g^{\prime a'b'} \Delta_{b'}$
and $\sum_{a''b''} \Delta_{a''} g^{\prime \prime a'' b''} \Delta_{b''}$
the Casimir elements for $g$ and $g'$.

\roster
\item "n=3" 
$$
(Y'_3 \otimes Y''_3)(\gamma'_1\otimes \gamma''_1 \otimes\gamma'_2\otimes 
\gamma''_2
\otimes\gamma'_3\otimes \gamma''_3) = 
Y'_3(\gamma'_1\otimes \gamma'_2 \otimes\gamma'_3)  Y'_3(\gamma''_1\otimes 
\gamma''_2
\otimes\gamma''_3) 
\tag 4.13
$$
\item "n=4" 
$$
\multline
(Y'_4 \otimes Y''_4)(\gamma'_1\otimes \gamma''_1 \otimes \dots 
\otimes\gamma'_4\otimes
\gamma''_4) = \\
Y'_4(\gamma'_1\otimes \dots \otimes\gamma'_4)  
\sum_{a'',b''}Y''_3(\gamma''_1\otimes
\gamma''_2 \otimes\Delta_{a''})  g^{\prime \prime a''b''} Y''_3 
(\Delta''_{b''}\otimes
\gamma''_3 \otimes \gamma''_4) +\\
\sum_{a',b'}Y'_3(\gamma'_1\otimes
\gamma'_2 \otimes\Delta_{a'})  g^{\prime a'b'} Y'_3 (\Delta'_{b'}\otimes
\gamma'_3 \otimes \gamma'_4)Y''_4(\gamma''_1\otimes \dots \otimes\gamma''_4)
\endmultline
\tag 4.14
$$

\item "n=5" 
$$
\multline
(Y'_5 \otimes Y''_5)(\gamma'_1\otimes \gamma''_1 \otimes \dots 
\otimes\gamma'_5\otimes
\gamma''_5) = \\
Y'_5(\gamma'_1\otimes \dots \otimes\gamma'_5)  
\sum_{a'',b'',c'',d''}Y''_3(\gamma''_1\otimes
\gamma''_2 \otimes\Delta_{a''})  g^{\prime \prime a''b''} Y''_3 
(\Delta''_{b''}\otimes
\gamma''_3 \otimes \Delta''_{c''})g^{\prime \prime c''d''} Y''_3 
(\Delta''_{d''}\otimes
\gamma''_4\otimes \gamma''_5)\\
- \sum_{l \in \{1,2,3,4\}} \sum \Sb a',b'\\a'',b'' \endSb
Y'_4(\bigotimes_{i \in  \{1,2,3,4\}\setminus \{l\}} \gamma'_i \otimes 
\Delta'_{a'}) g^{\prime
a'b'}
Y'_3(\Delta'_{b'} \otimes \gamma'_l \otimes \gamma'_5)\\
\times Y''_4(\bigotimes_{i \in  \{1,2,3,4\}\setminus \{l\}} \gamma''_i \otimes 
\Delta''_{a''}) g^{\prime \prime a''b''}
Y''_3(\Delta_{b''} \otimes \gamma''_l \otimes \gamma''_5)\\
+\sum_{\{1,2\} \subseteq I \subset \{1,2,3,4\}} \sum_{a',b'} 
Y'_{|I|+1}(\bigotimes_{i \in  I} \gamma'_i \otimes \Delta'_{a'}) g^{\prime a'b'}
Y'_{6-|I|} (\Delta'_{b'} \bigotimes_{j\in \{1,2,3,4\}\setminus I} \gamma'_j
 \otimes \gamma'_5)\\
\times \sum_{\{1,2\} \subseteq J \subset \{1,2,3,4\}} \sum_{a'',b''} 
Y''_{|J|+1}(\bigotimes_{i \in  J} \gamma''_i \otimes \Delta''_{a''}) 
g^{\prime a''b''}
Y''_{6-|I|} (\Delta''_{b''} \bigotimes_{j\in \{1,2,3,4\}\setminus J} \gamma''_j
 \otimes \gamma''_5)\\
+   \sum_{a',b',c',d'}Y'_3(\gamma'_1\otimes \gamma'_2 \otimes\Delta_{a'})  
g^{\prime a'b'}
Y'_3 (\Delta'_{b'}\otimes \gamma'_3 \otimes \Delta'_{c'})g^{\prime c'd'} Y'_3
(\Delta'_{d'}\otimes \gamma'_4\otimes \gamma'_5)Y''_5(\gamma''_1\otimes \dots
\otimes\gamma''_5)
\endmultline
\tag 4.15
$$
\endroster

\med

{\bf 4.4.2 Higher order multiplications}

\sm

By applying Corollary 4.2.5, 
using the notation $(\gamma_1, \ldots, \gamma_n)$
for $\circ_n(\gamma_1 \otimes \dots \otimes \gamma_n)$, we find:

\roster
\item "n=2"
$$ 
(\gamma'_1 \otimes \gamma''_1, \gamma'_2 \otimes \gamma''_2)
= (\gamma'_1 , \gamma'_2) \otimes (\gamma''_1 , \gamma''_2)
\tag 4.16
$$

\item "n=3"
$$ 
\multline
(\gamma'_1 \otimes \gamma''_1, \gamma'_2 \otimes \gamma''_2, 
\gamma'_3 \otimes \gamma''_3)= \\
(\gamma'_1 , \gamma'_2, \gamma'_3) \otimes ((\gamma''_1 , \gamma''_2), 
\gamma''_3) 
+ ((\gamma'_1 , \gamma'_2), \gamma'_3)\otimes (\gamma''_1, \gamma''_2, 
\gamma''_3)
\endmultline
\tag 4.17
$$

\item "n=4"
$$
\multline
(\gamma'_1 \otimes \gamma''_1 , \dots, \gamma'_4 \otimes \gamma''_4)= \\
(\gamma'_1 , \dots \gamma'_4) \otimes 
(((\gamma''_1,\gamma''_2),\gamma''_3),\gamma''_4)
+ (((\gamma'_1,\gamma'_2),\gamma'_3),\gamma'_4) \otimes (\gamma''_1 , \dots 
\gamma''_4)\\
- \sum_{\{i,j,k\} \amalg \{l\}= \{1,2,3,4\}} 
((\gamma'_i , \gamma'_j, \gamma'_k),\gamma'_l) \otimes ((\gamma''_i , 
\gamma''_j,
\gamma''_k),\gamma''_l)\\
+\sum _{\{1,2\} \subseteq I \subset \{1,2,3,4\}} 
((\gamma'_I),\gamma'_{\{1,2,3,4\}\setminus I})
\otimes \sum _{\{1,2\} \subseteq J \subset \{1,2,3,4\}}
((\gamma''_J),\gamma''_{\{1,2,3,4\}\setminus J}),
\endmultline
\tag 4.18
$$
where in the last expression we have used the abbreviation $(\gamma_I)$ 
to denote
$\circ_{|I|}(\otimes_{i\in I} \gamma_i)$.
\endroster

\newpage

\centerline{\bf References:}

\bigskip

[B] K.\ Behrend. {\it Gromov-Witten invariants in algebraic geometry.}
Preprint 1996.

\sm

[CP] C\. De Concini and C\. Procesi. {\it Wonderful models of 
subspace arrangements.}
Selecta Mathematica, New Series, 1 (1995), 459--494. 
{\it Hyperplane arrangements and holonomy equations.} ibid., 494--536.

\smallskip

[D] B.\ Dubrovin. {\it Geometry of 2D topological field theories}.
In: Springer LNM 1620 (1996), 120-348.

\sm

[G] E\. Getzler. {Operads and moduli spaces of genus zero Riemann surfaces.}
In: The Moduli Space of Curves, ed. by
R\. Dijkgraaf, C\. Faber, G\. van der Geer, Progress in Math\.
vol\. 129, Birkh\"auser, 1995, 199--230.

\sm

[GK] E.\ Getzler and M.\ M.\ Kapranov. {\it Modular operads.} Preprint 1994.

\sm

[K] M\. Kontsevich. {\it Enumeration of rational curves via torus actions.}
In: The Moduli Space of Curves, ed. by
R\. Dijkgraaf, C\. Faber, G\. van der Geer, Progress in Math\.
vol\. 129, Birkh\"auser, 1995, 335--368.

\smallskip

[Ke] S. Keel. {\it Intersection theory of moduli spaces of stable
$n$--pointed curves of genus zero.} Trans. AMS, 330 (1992), 545--574.

\smallskip

[KM] M\. Kontsevich and Yu\. Manin. {\it Gromov--Witten classes, quantum
cohomology, and enumerative geometry.} Comm. Math. Phys.,
164:3 (1994), 525--562.

\smallskip

[KMK] M\. Kontsevich and Yu\. Manin (with Appendix by R\. Kaufmann).
{\it Quantum cohomology of a product.} Invent\. Math\., 
124 (1996), 313--339.

\smallskip

[KMZ] R.\ Kaufmann, Yu.\ Manin and D.\ Zagier. {\it Higher Weil-Petersson
volumes of moduli spaces of stable $n$-pointed curves.} Preprint 1996

\sm

[Yu] S\. Yuzvinsky. {\it Cohomology basis for the De Concini--Procesi models of
hyperplane arrangements and sums over trees.} Preprint 1996.

\enddocument